\documentclass[12pt]{article}
\usepackage{amssymb,amsmath}
\newtheorem{DE}{Definition}[section]
\newtheorem{THEO}{Theorem}[section]
\newtheorem{RE}{Remark}[section]
\newtheorem{LE}{Lemma}[section]

\begin{document}
\title{Solutions of the Anisotropic Porous Medium Equation in $
R^n$ under $L^1$-initial Value
\thanks{Part of the work was done while B. Song was visiting the University de Autonoma de Madrid.
He thanks Professor Juan L. Vazquez for his invitation and  the
Mathematical Department for its hospitality and financial support.
The research is supported by National 973-project and the
Trans-Century Training Programme Foundation for the Talents by the
Ministry of Education.}}

\author{ { Huaiyu Jian and Binheng Song}\\ {  Department of Mathematical Sciences, }\\
{  Tsinghua University, Beijing 100084, China }  }
\date{}
\maketitle

{\bf Running head:}  Anisotropic Porous Medium Equation with
$L^1$-initial Value

\vskip 1cm

{\bf Correspondence Author and Address:}

Huaiyu Jian

Department of Mathematical Sciences

Tsinghua University

Beijing 100084, People's Republic of China

\vskip 0.4cm

{\bf E-mail:} hjian@math.tsinghua.edu.cn

{\bf Telephone:} 86-10-62772864

{\bf Fax:} 86-10-62781785

\newpage

 \centerline{\bf Abstract}

 Consider the anisotropic porous medium equation,
$u_t=\sum\limits_{i=1}^n(u^{m_i})_{x_ix_i},$ where $m_i>0,\
(i=1,2,\cdots,n)$ satisfying $\min\limits_{1\le i\le n}\{m_i\}\le
1,$ $\sum\limits_{i=1}^nm_i>n-2,$ and $\max\limits_{1\le i\le
n}\{m_i\}\le \frac{1}{n}(2+\sum\limits_{i=1}^nm_i).$ Assuming that
the initial data belong only to $L^1(\Re^n)$, we establish the
existence and uniqueness of the solution for the Cauchy problem in
the space, $C([0,\infty), L^1(\Re^n))\cap
C(\Re^n\times(0,\infty))\cap
L^\infty(\Re^n\times[\varepsilon,\infty)),$ where $\varepsilon>0$
may be arbitrary. We also show a comparison principle for such
solutions. Furthermore, we prove that the solution converges to
zero in the space $L^\infty(\Re^n)$ as the time goes to infinity.

\vskip 0.5cm {\bf Keywords:} Anisotropic diffusion, Degenerate
parabolic equation, Comparison principle, Large time behavior.
 \vskip 0.5cm
 {\bf AMS subject classifications:} 35K55, 35K65.

\section*{\centerline { 1. INTRODUCTION}}

\setcounter{section}{1}\setcounter{equation}{0}

 Consider the anisotropic porous
medium equation
\begin{equation}\displaystyle \label{1.1}u_t=\sum\limits_{i=1}^n(u^{m_i})_{x_ix_i}\quad \mbox{in }
\Re^n\times(0,\infty),\end{equation} where $m_i$ are positive
constants. As it is well-known, there have been a lot of works
dealing with the case of all $m_i$'s in (\ref{1.1}) being the same
positive constant, i.e., the case of porous medium equation (PME).
See, for example, the survey paper [1] for PME, the monograph [2]
for its generalization and the references there.

However, there are few papers on the general case of equation
(\ref{1.1}), although it has strong physical backgrounds. In fact,
it comes directly from water moves in anisotropic media. If the
conductivities of the media are different in different directions,
the constants $m_i$ in (\ref{1.1}) must be different from each
other. See [3] for details.

In papers [4, 5], the first author started studying the existence
and uniqueness for the Cauchy problem of equation (\ref{1.1}),
provided that the initial data are bounded and continuous in
$\Re^n.$ He also studied the continuous modulus of solutions to
(\ref{1.1}) in [6] (also see Lemma 3.1 in [4]). In [7], the
authors established the existence of fundamental solutions for the
Cauchy problem of equation (1.1).

In this paper, we will study the existence, uniqueness, comparison
principle and large time behavior of the solution of the Cauchy
problem for (\ref{1.1}) with $L^1$-initial data. For this purpose,
we want to consider the equation
\begin{equation}\displaystyle\label{1.2}V_{\tau}=\sum\limits_{i=1}^n\left[(V^{m_i})_{y_iy_i}
+\alpha_i(y_iV)_{y_i}\right]\mbox{ in }\Re^n\times(0,\infty)\end{equation}
and its stable equation
\begin{equation}\displaystyle \label{1.3}-\sum\limits_{i=1}^n\left[(f^{m_i})_{y_iy_i}
+\alpha_i(y_if)_{y_i}\right]=0\mbox{ in }\Re^n,\end{equation}
where $\alpha_i$ are defined by
\begin{equation}\displaystyle \label{1.4} \alpha_i=\frac{\overline m-m_i}{2}
+\frac{1}{n},\ (i=1,2,\cdots,n),\mbox{ and }\overline
m=\sum\limits_{i=1}^n\frac{m_i}{n}.\end{equation} As we see,
equation (1.2) is equivalent to (1.1), up to a scaling
transformation in spatial and time variables. See Lemma
\ref{lemma2} below.

For the coming needs, define
\begin{equation}\displaystyle \label{1.5}
Q=\Re^n\times(0,\infty),\quad\beta=\overline m-\frac{n-2}{n}
\end{equation}
and throughout this paper, we assume
\begin{equation}\displaystyle \label{1.6}
\min\limits_{1\le i\le n}\left\{m_i\right\}\le 1,\ \beta>0,\
m_i>0,\mbox{ and }\alpha_i>0\ (i=1,2,\cdots,n),\end{equation}
which is equivalent to
$$\label{1.6'}\begin{array}{l}\displaystyle
m_i>0\ (i=1,2,\cdots,n),\ \sum\limits_{i=1}^nm_i>n-2,\\[0.5cm]\displaystyle
\displaystyle \min\limits_{1\le i\le n}\left\{m_i\right\}\le
1,\mbox{ and }\max\limits_{1\le i\le n}\left\{m_i\right\}<
\frac{1}{n}(2+\sum\limits_{i=1}^nm_i).\end{array}\eqno(1.~6')$$
Also, we use $D'(X)$ to denote the set of all distributions
(generalized functions) in $X.$ The word ``{\sl respectively}" is
always shorten by ``{\sl resp}".

\begin{DE}
A $L^1_{loc}$-function $u(resp,V,f)$ is called a solution of
equation (\ref{1.1}) (resp, (1.2), (1.3)) if $u$ and $u^{m_i}$ are
in $D'(Q)$ (resp, $V$ and $V^{m_i}$ in $D'(Q)$, $f$ and $f^{m_i}$
in $D'(\Re^n)$) such that (1.1)    (resp, (1.2)), (1.3)) is
satisfied in the sense of distributions.
\end{DE}

\begin{DE}
An $L^1$-regular solution, $u$, of (\ref{1.1}) (resp, (\ref{1.2}))
with initial value $u(x,0)=u_0(x)$ in $\Re^n $ is a solution of (\ref{1.1})
(resp, (\ref{1.2})) satisfying
\begin{equation}\displaystyle \label{1.7}
u\in C([0,\infty),L^1(\Re^n))\bigcap C(Q)\bigcap
L^\infty(\Re^n\times[\varepsilon,\infty))\end{equation} for each
$\varepsilon>0,$ and
\begin{equation}\displaystyle \label{1.8}
\displaystyle\int\limits_{\Re^n}u(x,t)dx=\int\limits_{\Re^n}u_0(x)dx,\
\forall t>0.\end{equation}\end{DE}

\begin{DE}
If we replace ``$=$" by ``$\ge$" (resp, ``$\le$") in (\ref{1.1}),
(\ref{1.2}), (\ref{1.3}) and (\ref{1.8}), then we obtain the
definitions of the super-solution (resp, sub-solution).
\end{DE}

For simplicity, we will consider only the nonnegative solutions.
The main results of this paper may be stated as the following
theorems for which we always assume $(1.~6').$

\begin{THEO}\label{theo1}
Suppose that $\underline u$ and $\overline u$ are nonnegative
$L^1$-regular sub-solution and super-solution, respectively, of
(\ref{1.1}) (resp, of (\ref{1.2})). If $\underline
u(x,0)\le\overline u(x,0)$ in $\Re^n,$ then $\underline
u\le\overline u$ in $Q.$\end{THEO}
\begin{THEO}\label{theo2}
For any nonnegative function $u_0\in L^1(\Re^n),$ there exists a
unique nonnegative $L^1$-regular solution $u$ of (\ref{1.1})
(resp, (\ref{1.2})) with initial value $u(x,0)=u_0(x)$ in $\Re^n.$
Furthermore, $u$ satisfies the following decay estimate:
$$\displaystyle\left\|u(x,t)\right\|_{L^\infty(\Re^n)}\le Ct^{-\frac{1}{\beta}}
\left\|u_0\right\|_{L^1(\Re^n)}^{\frac{2}{\beta}}, \ \forall
t>0,$$
where $C$ is a constant depending only on $m_i\
(i=1,2,\cdots,n)$ and $n$.\end{THEO}

Note that the decay estimate says that the solution of (\ref{1.1})
converges to zero uniformly for $x\in\Re^n$ as the time goes to infinity.

We would like to mention that the Cauchy problem of PME with
$L^1$-initial data was studied in [8, 9]. But our methods to prove
Theorems \ref{theo1} and \ref{theo2} are completely different from
those in [8, 9]. In fact, the proof of Theorem \ref{theo1} is
based on a max-min method, which will be developed in Section 2.
As a key step for the proof of Theorem \ref{theo2}, we have to
discover a specific scaling technique for each direction in order
to overcome the difficulty from the anisotropic phenomenon so that
we can construct suitable super-solution, which will be developed
in Section 3. While in Section 4, we will combine approximation
arguments and the results in Section 3 to prove Theorem
\ref{theo2}.

 \section*{\centerline{2. A COMPARISON PRINCIPLE }}
\setcounter{section}{2} \setcounter{equation}{0}
 Theorem
\ref{theo1} is a comparison principle for $L^1$-regular solutions
of (\ref{1.1}). In this section, we will give its proof by a few
lemmas. The following result was proved in [4, 5].

\begin{LE}\label{lemma1}

There exists a unique bounded and continuous solution of the
Cauchy problem to (\ref{1.1}) if the initial value is in
$C(\Re^n)\bigcap L^\infty(\Re^n).$  The solution is approximated
by a equi-continuous and uniformly bounded sequence of classical
solutions in the norm of the space
$C_{loc}(\Re^n\times[0,\infty))$ and thus satisfies a comparison
principle .
\end{LE}

\begin{LE}\label{lemma2}
Suppose that the functions $u(x,t)$ and $V(y,\tau)$ are related by
\begin{equation}\displaystyle \label{2.1}
V(y,\tau)=h(t)u(x,t),\ \tau=\ln h(t),\ y_i=x_ih^{-\alpha_i}(t),
(i=1,2\cdots,n),\end{equation} where $\alpha_i$ and $\beta$ are
the same as in (\ref{1.4})--(\ref{1.6}), and
\begin{equation}\displaystyle \label{2.2}h(t)=(1+\beta
t)^{\frac{1}{\beta}}\end{equation} If $u(x,t)$ is a solution
(resp, super-solution, sub-solution) of (\ref{1.1}), then
$V(y,\tau)$ is a solution (resp, super-solution, sub-solution) of
(\ref{1.2}), and vice versa.
\end{LE}

\noindent{\bf Proof.} It is from direct computations. we omit the
details.
\begin{LE}\label{lemma3}
If $u$ and $w$ are nonnegative sub-solutions (resp,
super-solutions) of (\ref{1.1}) satisfying
\begin{equation}\displaystyle \label{2.3}
u,w,u_t,w_t,(u^{m_i})_{x_ix_i},(w^{m_i})_{x_ix_i}\in L^1_{loc}(Q)\
(i=1,2,\cdots,n),\end{equation} then $\max\{u,w\}$ (resp,
$\min\{u,w\}$) is a sub-solution (resp, super-solution) of
(\ref{1.1}).
\end{LE}

\noindent{\bf Proof.} First suppose that $u$ and $w$ are
sub-solution of (\ref{1.1}). Note that $$\displaystyle
(\left|g\right|)_t=\mbox{sign}(g)g_t$$ in the sense of
distributions if $g$ and $g_t$ are in $L^1_{loc}.$ Furthermore, by
Kato's inequality, we have, in the same sense, that
$$\displaystyle (\left|g\right|)_{x_ix_i}\ge\mbox{sign}(g)g_{x_ix_i}$$
if $g$ and $g_{x_ix_i}$ are in $L^1_{loc}(Q)$. It follows from
these two facts and (\ref{2.3}) that the following computations
hold true in the sense of distributions:

\[\begin{array}{l}
\displaystyle\frac{\partial}{\partial t}\max\left\{u,w\right\}-\sum\limits_{i=1}^n
\frac{\partial^2}{\partial x_i^2}\left(\max\left\{u,w\right\}\right)^{m_i}\\[0.5cm]\displaystyle
\displaystyle=\frac{\partial}{\partial t}\max\left\{u,w\right\}-\sum\limits_{i=1}^n
\frac{\partial^2}{\partial x_i^2}\max\left\{u^{m_i},w^{m_i}\right\}\\[0.5cm]\displaystyle
\displaystyle=\left(\frac{u+w}{2}+\frac{\left|u-w\right|}{2}\right)_t-\sum\limits_{i=1}^n
\left(\frac{u^{m_i}+w^{m_i}}{2}+\frac{\left|u^{m_i}-w^{m_i}\right|}{2}\right)_{x_ix_i}\\[0.5cm]\displaystyle
\displaystyle\le\left(\frac{u+w}{2}\right)_t+\mbox{sign}\left(u-w\right)\left(\frac{u-w}{2}\right)_t\\[0.5cm]\displaystyle
\displaystyle\quad-\sum\limits_{i=1}^n
\left(\frac{u^{m_i}+w^{m_i}}{2}\right)_{x_ix_i}
-\sum\limits_{i=1}^n
\mbox{sign}\left(u-w\right)\left(\frac{u^{m_i}-w^{m_i}}{2}\right)_{x_ix_i}\\[0.5cm]\displaystyle
=\left\{\begin{array}{ll}u_t-\sum\limits_{i=1}^n\left(u^{m_i}\right)_{x_ix_i},&\mbox{ if }u\ge w\\[0.5cm]\displaystyle
w_t-\sum\limits_{i=1}^n\left(w^{m_i}\right)_{x_ix_i},&\mbox{ if }u<w
\end{array}\right.\\[0.5cm]\displaystyle
\le 0.
\end{array}\]
This shows the maximum of two sub-solutions is also a
sub-solution. Similarly, the minimum of two super-solutions is a
super-solution.

\begin{LE}\label{lemma4}
If $u$ and $w$ are nonnegative solutions of (\ref{1.1}) such that
\begin{equation}\displaystyle \label{2.4}
u,w,u^{m_i},w^{m_i}\in L^1_{loc}(Q)
\end{equation}
and
\begin{equation}\displaystyle \label{2.5}
\lim\limits_{k\rightarrow\infty}(u_k,w_k,u_k^{m_i},w_k^{m_i})=(u,w,u^{m_i},w^{m_i})\in
\left[L^1_{loc}(Q)\right]^4
\end{equation}
for $i=1,2,\cdots,m$ and some nonnegative function sequence
$\{u_k\}$ and $\{w_k\}$ satisfying (\ref{2.3}), then $\max\{u,w\}$
(resp, $\min\{u,w\}$) is a sub-solution (resp, super-solution) of
(\ref{1.1}).
\end{LE}

\noindent{\bf Proof.} We want only to consider the maximum case.
Let
$$U_k=\max\{u_k,w_k\}, U=\max\{u, w \},\ \ $$ Since $u_k$ and $w_k$ satisfy
(\ref{2.3}), we have, by Lemma \ref{lemma3}, that
\[\displaystyle
\int_0^\infty\int\limits_{\Re^n}\left(U_k\varphi_t+\sum\limits_{i=1}^n
U_k^{m_i}\varphi_{x_ix_i}\right)dxdt\ge 0
\]
for all nonnegative function $\varphi\in C_0^\infty(Q)$
and all $k=1,2,\cdots$. Using (\ref{2.4}) and (\ref{2.5}), and passing the limit, we obtain that
\[\displaystyle
\int_0^\infty\int\limits_{\Re^n}\left(U\varphi_t+\sum\limits_{i=1}^n
U^{m_i}\varphi_{x_ix_i}\right)dxdt\ge 0,\ \forall\varphi\in
C_0^\infty(Q),\varphi\ge 0.
\]
This shows that $U$ is a sub-solution of (\ref{1.1}).

\underline{\bf Proof of Theorem \ref{theo1}:} Owing to Lemma
\ref{lemma2}, we want only to prove Theorem \ref{theo1} for
equation (\ref{1.1}). We will complete the proof by three steps.

\underline{\bf Step 1.} Suppose that $\overline u$ and $u$ are,
respectively, nonnegative $L^1$-regular super-solution and
$L^1$-regular solution of (\ref{1.1}) with $\overline u(x,0)\ge
u(x,0)$ in $\Re^n.$ We will prove that
\begin{equation}\displaystyle \label{2.6}\overline u\ge u\mbox{ in }Q.\end{equation}

Fix $\tau>0.$ Since $\overline u(\cdot,\tau)\in C(\Re^n)\bigcap
L^\infty(\Re^n)$, it follows from Lemma 2.1 that there exists a
unique solution $u^{(\tau)}\in C(Q)\bigcap L^\infty(Q)$ of
(\ref{1.1}) with the initial value $u^{(\tau)}(x,0)=\overline
u(x,\tau)\in \Re^n,$ which is approximated by a sequence of
classical solutions in the norm of
$C_{loc}(\Re^n\times[0,\infty))$. The same conclusion holds true
for the solution $V(x,t)=u(x,t+\tau)$ since $u(\cdot,\tau)\in
C(\Re^n)\bigcap L^\infty(\Re^n)$. Hence Lemma \ref{lemma4} implies
that $\min\{u^{(\tau)}(x,t),u(x,t+\tau)\}$ is a bounded and
continuous super-solution of (\ref{1.1}). This implies
\begin{equation}\displaystyle \label{2.7}\displaystyle
\int\limits_{\Re^n}\min\left\{u^{(\tau)}(x,t),u(x,t+\tau)\right\}dx\ge
\int\limits_{\Re^n}\min\left\{\overline
u(x,\tau),u(x,\tau)\right\}dx,  \ \forall t>0.\end{equation}

On the other hand, $\overline u(\cdot,\tau+\cdot)$ is a bounded
and continuous super-solution of (\ref{1.1})
 with the same initial value in
$C(\Re^n)\bigcap L^\infty(\Re^n)$ as $u^{(\tau)}$. Hence, by the
comparison principle of Lemma 2.1 we have
$$\displaystyle \overline u(x,t+\tau)\ge u^{(\tau)}(x,t), \forall x\in\Re^n,\ \forall t>0,$$
 which, together with (\ref{2.7}), implies
\[\displaystyle
\int\limits_{\Re^n}\min\left\{\overline
u(x,t+\tau),u(x,t+\tau)\right\}dx\ge
\int\limits_{\Re^n}\min\left\{\overline
u(x,\tau),u(x,\tau)\right\}dx, \ \forall t>0.\] Now letting
$\tau\rightarrow 0$, using the fact of $u$ and $\overline u$ in
$C([0,\infty),L^1(\Re^n))$ and observing that $u$ satisfies
(\ref{1.8}) we see that
\[\begin{array}{l}
\displaystyle\int\limits_{\Re^n}\min\left\{\overline u(x,t),u(x,t)\right\}dx\\[0.5cm]\displaystyle
\displaystyle\ge\int\limits_{\Re^n}\min\left\{\overline u(x,0),u(x,0)\right\}dx\\[0.5cm]\displaystyle
\displaystyle=\int\limits_{\Re^n}u(x,0)dx =
\int\limits_{\Re^n}u(x,t)dx, \forall t>0.\end{array}\] This
implies the desired result (\ref{2.6}).

\underline{\bf Step 2.} Suppose that $\underline u$ and $u$,
respectively, nonnegative $L^1$-regular sub-solution and
$L^1$-regular solution of (\ref{1.1}) with $\underline u(x,0)\le
u(x,0)$ in $\Re ^n$. In this case, we fix $\tau>0$ and consider
the unique solution $u^{(\tau)}\in C(Q)\bigcap L^\infty(Q)$ with
$u^{(\tau)}(x,0)=\underline u(x,\tau )$. Then the same arguments
as in {\sl Step 1} (replaced ``min" by ``max") yields $\underline
u\le u$ in $Q$.

\underline{\bf Step 3.}  Suppose that $\underline u$ and
$\overline u$, respectively, nonnegative $L^1$-regular
sub-solution and $L^1$-regular super-solution. Then fix $\tau>0$
and let $u$ be the unique $L^1$-regular solution of (\ref{1.1})
with initial value $u(x,t)=\underline u(x,\tau)$ at $t=\tau$. Note
that the existence and the uniqueness of such a solution  follow
from Theorem \ref{theo2}, whose proof is independent of the {\sl
Step 3} (but depends on the {\sl Step 1} and {\sl Step 2}).
Combining the results of Steps 1 and 2, we have
$$\displaystyle \underline u(x,t)\le u(x,t)\le\overline u(x,t),\ \forall x\in\Re^n,\ \forall t>\tau.$$
This proves Theorem \ref{theo1} since $\tau$ is arbitrary.

\section*{\centerline{3.  SOLUTIONS WITH INITIAL DATA in {$C_0^{\infty}(\Re ^n)$}}}
\setcounter{section}{3}\setcounter{equation}{0}

 In the proof of
Theorem \ref{theo2}, we will approximate the initial value of the
solution to (\ref{1.1}) by a sequence of functions in
$C_0^\infty(\Re^n)$. Thus, the first step is to study the
solutions whose initial data are compactly supported in $\Re^n$.
In this section, we will construct a super-solution of (\ref{1.1})
by a scaling technique and then use the super-solution to study
such solutions.

\vskip 0.3cm \noindent {\bf Lemma 3.1}  {\em  Suppose that $u_0\in
C_0^\infty(\Re^n)$ and let $Q_T=\Re^n\times[0,T)$. Then for any
$T>0,$ there exists a positive function $\overline u\in
C(Q_T)\bigcap L^\infty(Q_T)\bigcap C([0,T),L^1(\Re^n))$ such that
$\overline u$ is a super-solution of (\ref{1.1}) in the domain
$Q_T$ and $u_0(x)\le \overline u(x,0)$ in $\Re^n.$ }

\noindent{\bf Proof.} Since $\beta>0$ (recall
(\ref{1.4})--($1.6'$)), $\overline m>\frac{n-2}{n}$ and thus
$\frac{1-m_i}{2}<\alpha_i$ for each $i$. Consequently, we can
choose $\theta_i>2$ and $\alpha>0$ such that
\begin{equation}\displaystyle \label{3.1}
\frac{1-m_i}{2}<\frac{1}{\alpha\theta_i}<\alpha_i,\
(i=1,2,\cdots,n).
\end{equation}
Hence
\begin{equation}\displaystyle \label{3.2}
\sum\limits_{i=1}^n\alpha_i=1,\mbox{ and
}-m_i\alpha-\frac{2}{\theta_i}=-m_i\alpha-2+\frac{2\theta_i-2}{\theta_i}<-\alpha.
\end{equation}
Let $\mu_i=\frac{1-m_i}{2},$
\[R_0=\max\left\{1,\left[\frac{n\max\limits_{1\le i\le n}
\left\{m_i\alpha(m_i\alpha+1)\theta_i^2\right\}}{\min\limits_{1\le i\le n}
\left\{\alpha_i\alpha\theta_i\right\}-1}\right]^{\frac{1}{2\alpha\max\limits_{1\le i\le n}
\left\{(\theta_i\alpha)^{-1}-\mu_i\right\}}}\right\}\]
and $$\displaystyle \Omega_{R_0}=\left\{y\in\Re^n:\sum\limits_{i=1}^n
\left|y_i\right|^{\theta_i}>R_0\right\}.$$
Then the function
$$\displaystyle f(y)=\left(\sum\limits_{i=1}^n\left|y_i\right|^{\theta_i}
\right)^{-\alpha}$$ belongs to $C^2(\Omega_{R_0})\bigcap
L^1(\Omega_{R_0})$ and satisfies $$\displaystyle
\sum\limits_{i=1}^n\left[\left(f^{m_i}\right)_{y_i y_i}
+\alpha_i\left(y_if\right)_{y_i}\right]\le 0 \mbox{ in }
\Omega_{R_0}$$ by a direct computation using (\ref{3.2}) and
(\ref{1.6}).   The last inequality tells us that $f$ is a
super-solution of (\ref{1.3}) in the domain $\Omega_{R_0}.$ This
implies that for any constant $\lambda,$ the function
\begin{equation}\displaystyle \label{3.3}
f^{(\lambda)}(y)=\lambda^{\frac{2}{n\beta}}\left(\sum\limits_{i=1}^n
\left|y_i\right|^{\theta_i}\lambda^{\frac{1-m_i}{n\beta}\theta_i}\right)^{-\alpha}
\end{equation} is also a super-solution of (\ref{1.3}) in the domain
$$\displaystyle \Omega^{(\lambda)}_{R_0}=\left\{y=(y_1,y_2,\cdots,y_n):y_i
=\lambda^{\frac{m_i-1}{n\beta}}x_i,x=(x_1,x_2,\cdots,x_n)\in \Omega_{R_0}\right\}.$$
Now for any positive constants $T,\ C_0$ and $A,$ choose a $\lambda$
such that
$$\displaystyle \lambda\ge\max\left\{\left[C_0\left(1+\beta T\right)^{\frac{1}{\beta}}R_0^{\alpha}
\right]^{\frac{n\beta}{2}},
\left(C_0A^\alpha\right)^{\frac{n\beta}{2}}\max\limits_{1\le i\le
n} \left(\alpha^{-1}-\mu_i\theta_i\right)^{-1}\right\}$$ and let
$$\displaystyle Q(C_0,T)=\left\{(x,t)\in Q_T:(1+\beta t)^{-\frac{1}{\beta}}
f^{(\lambda)}\left(x_1(1+\beta
t)^{-\frac{\alpha_1}{\beta}},\cdots,x_n(1+\beta
t)^{-\frac{\alpha_n}{\beta}} \right)\le C_0 \right\}.$$
 We may
verify that
$$\displaystyle \left\{y=(y_1,y_2,\cdots,y_n):y_i=x_i(1+\beta t)^{-\frac{\alpha_i}{\beta}},\ (x,t)\in Q(C_0,T)
 \right\}\subset
 \Omega^{(\lambda)}_{R_0}$$ and the function
\begin{equation}\displaystyle \label{3.4}
\overline u(x,t)=\left\{\begin{array}{ll}
\!(\!1\!+\!\beta t)\!^{-\frac{1}{\beta}}
\!f\!^{(\lambda)}\!\left(x_1\!(1\!+\!\beta t)^{-\frac{\alpha_1}{\beta}},\cdots,
x_n\!(\! 1\!+\!\beta\! t\!)^{-\frac{\alpha_n}{\beta}}\!\right)
\!,& \!\!(\!x,\!t\!)\!\in\! Q\!(C_0,\!T)\!\\[0.5cm]\displaystyle
C_0, &\!\! (x,\!t)\!\!\in\!\! Q_T\!\!\setminus\!\!
Q\!(C_0,T\!)\end{array}\right.\end{equation} is a super-solution
of (\ref{1.1}) in $Q_T$. Obviously, $\overline u\in C(Q_T)\bigcap
L^\infty(Q_T)\bigcap C([0,T),L^1(\Re^n)).$ Moreover, we easily see
that
$$\displaystyle (1+\beta t)^{-\frac{1}{\beta}}
f^{(\lambda)}\left(x_1(1+\beta t)^{-\frac{\alpha_1}{\beta}},\cdots,
x_n( 1+\beta t)^{-\frac{\alpha_n}{\beta}}\right)=C_0$$
for all $t\in[0,T)$ and all the $x=(x_1,x_2,\cdots,x_n)$ satisfying
$\left|x_i\right|\le A^{\frac{1}{\theta_i}},\ i=1,2,\cdots,n.$
Hence, $$\displaystyle \overline u(x,t)=C_0,\ \forall(x,t)\in
\left(\prod\limits_{i=1}^n[-A^{\frac{1}{\theta_i}},A^{\frac{1}{\theta_i}}]\right)\times[0,T].$$
Consequently, choose $A$ and $C_0$ large enough such that
$$\displaystyle \prod\limits_{i=1}^n[-A^{\frac{1}{\theta_i}},A^{\frac{1}{\theta_i}}]
\supset\mbox{supp}u_0\mbox{ and }C_0\ge\max_{x\in \Re ^n}
u_0(x),$$ we have proved the lemma.

\begin{RE}
Using (\ref{3.4}) and Lemma 2.2, we see that the function
$$\displaystyle \overline V(y,\tau)=\left\{\begin{array}{ll}
f^{(\lambda)}(y),&(y,\tau)\in P=\left\{(y,\tau):f^{(\lambda)}(y)
\le C_0e^\tau,0\le\tau\le\frac{\ln(1+\beta T)}{\beta}\right\}\\[0.5cm]\displaystyle
C_0e^\tau,&\mbox{otherwise}\end{array}\right.$$ is a
super-solution of (\ref{1.2}).
\end{RE}

\vskip 0.3cm \noindent {\bf Lemma 3.2} {\em Suppose that $u_0\in
C_0^\infty(\Re^n),\ u_0\ge 0$ in $\Re^n.$ Let $u$ be the unique
bounded and continuous solution with $u(x,0)=u_0(x)$ in $\Re^n.$
Then one has
\begin{equation}\displaystyle \label{3.5}
u\in C([0,\infty),L^1(\Re^n)),
\end{equation}
\begin{equation}\displaystyle \label{3.6}
\int\limits_{\Re^n}u(x,t)dx=\int\limits_{\Re^n}u_0(x)dx, \ \forall
t>0
\end{equation}
and
\begin{equation}\displaystyle \label{3.7}
\left\|u(\cdot,t)\right\|_{L^\infty(\Re^n)}\le
Ct^{-\frac{1}{\beta}}
\left\|u_0\right\|_{L^1(\Re^n)}^{\frac{2}{\beta}}, \ \forall t>0,
\end{equation}
where $C$ is a constant depending only on $n$ and $m_i\
(i=1,2,\cdot,n).$ }

 {\bf Proof.} Fix $T>0.$ By Lemma 3.1 and Lemma 2.1
  we see that
\begin{equation}\displaystyle \label{3.8}
0\le u(x,t)\le\overline u(x,t)\mbox{ in }Q_T,
\end{equation}
where $\overline u$ is the same as in (\ref{3.4}). Hence (\ref{3.5})
comes from the continuity of $u$ and the structure of $\overline u.$

(\ref{3.7}) was proved in [6]. To show (\ref{3.6}), we use Lemma
2.1 to see  that $u$ can be approximated by an (equi-continuous
and uniformly bounded) sequence $\{u_k\}_{k=1}^\infty$ of
classical solutions of (\ref{1.1}) in the norm space of the space
$C_{loc}(\Re^n\times[0,\infty)).$ Multiplying (\ref{1.1}) for
$u_k$ by $\phi\in C_0^\infty(\Re^n),$ integrating the resulted
equation for $x$ over $\Re^n$ and for $t$ over $[0,T],$ and
passing the limit as $k\rightarrow\infty,$ we obtain that
\begin{equation}\displaystyle \label{3.9}
\int\limits_{\Re^n}u(x,T)\phi(x)dx-\int\limits_{\Re^n}u(x,0)\phi(x)dx
=\int_0^T\int\limits_{\Re^n}\sum\limits_{i=1}^nu^{m_i}\phi_{x_ix_i}dxdt.
\end{equation}
Choose a function $g\in C_0^\infty(-\infty,\infty)$ satisfying
$$\displaystyle 0\le g \le1,\ \left|g'\right|\le 4,\ \left|g''\right|\le 8\mbox{ in } R$$
and $$\displaystyle g(s)=\left\{\begin{array}{ll}
1,&\left|s\right|\le 1\\[0.5cm]\displaystyle0,&2\le\left|s\right|\le\infty.\end{array}\right.$$
Taking $\phi$ in (\ref{3.9}) as
$$\displaystyle \phi^{(k)}(x)=\prod\limits_{i=1}^ng(x_1k^{-\frac{1}{\theta_1}},
\cdots,x_nk^{-\frac{1}{\theta_n}}),\ k=1,2,\cdots$$
and observing that
$$\displaystyle \left|\int_0^T\int\limits_{\Re^n}\sum\limits_{i=1}^nu^{m_i}
\phi^{(k)}_{x_ix_i}dxdt\right|\le C(T,n,m_i)k^{-\min\limits_{1\le
i\le n}\{\theta_i^{-1}\}} \rightarrow 0$$ as $k\rightarrow\infty,$
where we have used (\ref{3.8}) and (\ref{3.4}), we finally find
that
$$\displaystyle \int\limits_{\Re^n}u(x,T)dx=\int\limits_{\Re^n}u(x,0)dx.$$
This proves (\ref{3.6}) since $T$ is arbitrary.

\vskip 0.3cm
 \noindent {\bf Lemma 3.3} {\em Suppose that $u(x,t)$
is defined as in Lemma 3.2 and $V(y,\tau)$ is related to $u(x,t)$
as in Lemma 2.2. Then there exists a function $F\in
C(\Re^n)\bigcap L^\infty(\Re^n)\bigcap L^1(\Re^n),$ depending only
on the upper bounds of $\left\|u_0\right\|_{L^\infty(\Re^n)}$ and
$\left\|u_0\right\|_{L^1(\Re^n)},$ such that
\begin{equation}\displaystyle \label{3.10}
0\le V(y,\tau)\le F(y),\ \forall y\in\Re^n,\
\forall\tau>0.\end{equation} }

{\bf Proof.} Owing (\ref{3.8}), (\ref{3.3}), (\ref{3.4}) and
(\ref{3.7}), we see that there is a constant $C_1,$ depending only
on $n,m_i,\left\|u_0\right\|_{L^\infty(\Re^n)}$ and
$\left\|u_0\right\|_{L^1(\Re^n)},$ such that $$\displaystyle 0\le
V(y,\tau)\le C_1,\ \forall (y,\tau)\in Q.$$

Since $V_0(y)=V(y,0)=u_0(y)\in C_0^\infty(\Re^n),$ choose a
$R_1\ge R_0$ such that
$\mbox{supp}V_0\subset\Re^n\setminus\Omega_{R_1}.$ Repeating the
arguments from (\ref{3.2}) to (\ref{3.3}), we find a $\lambda_1>0$
and a bounded and continuous super-solution $f^{(\lambda_1)}$ of
(\ref{1.3}) in the domain $\Omega^{\lambda_1}_{R_1}$ such that
$f^{(\lambda_1)}>0$ in $\Omega_{R_1}^{(\lambda_1)}$ and
$f^{(\lambda_1)}=C_1$ on $\partial\Omega_{R_1}^{(\lambda_1)}.$ Let
$$\displaystyle F(y)=\left\{\begin{array}{ll}
f^{(\lambda_1)}(y),&y\in\Omega_{R_1}^{(\lambda_1)}\\[0.5cm]\displaystyle
C_1,&y\in\Re^n\setminus\Omega_{R_1}^{(\lambda_1)}.\end{array}\right.$$
Then $F$ is a bounded and continuous super-solution of (\ref{1.3})
in $\Re^n$ (of (\ref{1.1}) in $\Re^n\times(0,\infty)$) satisfying
$F\ge V_0$ in $\Re^n.$  Consequently, the comparison principle
(Lemma 2.1) implies (\ref{3.10}).

\section*{\centerline{4. PROOF OF THEOREM  \ref{theo2}}}

\setcounter{section}{4}\setcounter{equation}{0} We want to prove
Theorem \ref{theo2} only for equation (\ref{1.1}) due to Lemma
2.2. Note that the uniqueness is direct from the Steps 1 and 2 in
the proof of Theorem \ref{theo1} and the proof of existence below
is independent of Step 3 there.

\underline{\bf Step 1.} Suppose the initial value $u_0\in L^\infty(\Re^n)$
is nonnegative such that $\mbox{supp}u_0$ is a bounded set in $\Re^n.$
We will prove that (\ref{1.1}) has a nonnegative $L^1$-regular solution $u$
satisfying $u(x,0)=u_0(x)$ in $\Re^n.$

Choose a sequence $u_{0k}\subset C^\infty_0(\Re^n)$ such that
\begin{equation}\displaystyle \label{4.1}\left\|u_{0k}\right\|\le C_2,\ u_{0k}\le u_{0k'}\le u_{0}\mbox{ in }\Re^n,
\ \mbox{supp}u_{0k}\subset\Omega,\ (k\le k',\ \
k,k'=1,2,\cdots)\end{equation} and
\begin{equation}\displaystyle \label{4.2}\lim\limits_{k\rightarrow\infty}u_{0k}=u_0
\mbox{ in }L^\infty(\Re^n)\bigcap L^1(\Re^n), \end{equation} where
constant $C_2$ and domain $\Omega\subset\Re^n$ are independent of
$k.$ For each $k$, Let $u_k$ be the solution of (\ref{1.1}) with
initial value $u_{0k}.$ As we have said, Lemma 2.1 imply that each
$u_k$ is a unique, continuous and bounded in $\Re^n,$ satisfying
\begin{equation}\displaystyle \label{4.3}\left\|u_{k}\right\|_{L^\infty(Q)}\le C_3,\ 0\le u_{k}\le u_{k'}\mbox{ in }Q,
 \ (k\le k', \ \ k,k'=1,2,\cdots)\end{equation}
for some constant $C_3$ depending on $C_2$ but independent of $k$
(also see (3.7)), and
\begin{equation}\label{4.4}\begin{array}{l}\displaystyle \int\limits_{\Re^n}u_k(x,t)dx
=\int\limits_{\Re^n}u_{0k}(x)dx,\ \forall t>0\\[0.5cm]\displaystyle\mbox{ and }
\lim\limits_{t\rightarrow 0}u_k(x,t)=u_{0k}(x) \mbox{ in
}L^1(\Re^n),\ (k=1,2,\cdots)\end{array}\end{equation} by
(\ref{3.5}) and (\ref{3.6}). Hence, $\{u_k\}$ is equi-continuous
in $\Re^n\times[\tau_0,\infty)$ for any $\tau_0>0$ by Lemma 3.1 in
[4] (or more generally, Theorem 5.1 in [6]).

Now fix $\tau_0>0.$ Let $V_k(y,\tau)$ be the solutions of
(\ref{1.2}) related to $u_k$ by Lemma 2.2. Then by Lemma 3.3 we
have a function
\begin{equation}\displaystyle \label{4.5}
F\in C(\Re^n)\bigcap L^\infty(\Re^n)\bigcap L^1(\Re^n)
\end{equation}
such that
\begin{equation}\displaystyle \label{4.6}0\le V_k(y,\tau)\le V_{k'}(y,\tau)\le F(y),\
\forall y\in\Re^n,\ \forall\tau>0,\ (k\le k', \ \
k,k'=1,2,\cdots).\end{equation} Therefore, we can choose a
subsequence of $\{V_k\}$ such that
\begin{equation}\displaystyle \label{4.7}\lim\limits_{k\rightarrow\infty}
V_k(y,\tau)=V(y,\tau),\ \forall
(y,\tau)\in\Re^n\times[\tau_0,\infty)\end{equation} and
\begin{equation}\displaystyle
\label{4.8}\lim\limits_{k\rightarrow\infty}
V_k(\cdot,\tau)=V(\cdot,\tau)\mbox{ in }L^1(\Re^n),\
\forall\tau\ge\tau_0\end{equation} for some function $V\in
C(\Re^n\times[\tau_0,\infty))$ satisfying
\begin{equation}\displaystyle \label{4.9}0\le V_k(y,\tau)\le V(y,\tau)\le F(y),\
\forall y\in\Re^n,\
\forall\tau>\tau_0,(k=1,2,\cdots).\end{equation} Obviously, $V$ is
solution of (1.2) satisfying
\begin{equation}\displaystyle \label{4.10}
V\in C(\Re^n\times(0,\infty))\bigcap L^\infty(\Re^n\times(\tau_0,\infty))
\end{equation} since $\tau_0$ is arbitrary.

For any $t_0>0$, choose $\tau_0\in(0,t_0).$ Then it follows from (\ref{4.9})
and the continuity of $V$ that
\begin{equation}\displaystyle \label{4.11}\lim\limits_{\tau\rightarrow t_0}
 V(\cdot,\tau)=V(\cdot,t_0)\mbox{ in }L^1(\Re^n).\end{equation}
It follows from Lemma 2.2, (\ref{4.4}) and the fact of
$\sum\limits_{i=1}^n\alpha_i=1$ that
\begin{equation}\displaystyle \label{4.12}\int\limits_{\Re^n}V_k(y,\tau)dy
=\int\limits_{\Re^n}u_{k}(x,t)dx=\int\limits_{\Re^n}u_{0k}(x)dx,\
k=1,2,\cdots,\ \forall t>0,\ \tau=\ln h(t),\end{equation} which,
together (\ref{4.2}), (\ref{4.8}) and the fact of
$\sum\limits_{i=1}^n\alpha_i=1$ again, implies
\begin{equation}\displaystyle \label{4.13}\int\limits_{\Re^n}u(x,t)dx
=\int\limits_{\Re^n}V(y,\tau)dy=\int\limits_{\Re^n}u_{0}(x)dx,\ \forall t>0,\ \tau=\ln h(t).\end{equation}
Note that (\ref{4.9}), (\ref{4.12}) and (\ref{4.13}) imply
\[\begin{array}{l}
\left\|u(\cdot,t)-u_0\right\|_{L^1(\Re^n)}\\[0.5cm]\displaystyle
\le\left\|u(\cdot,t)-u_k(\cdot,t)\right\|_{L^1(\Re^n)}+
\left\|u_k(\cdot,t)-u_{0k}\right\|_{L^1(\Re^n)}+\left\|u_{0k}-u_0\right\|_{L^1(\Re^n)}\\[0.5cm]\displaystyle
=\int\limits_{\Re^n}\left(u(x,t)-u_k(x,t)\right)dx+
\left\|u_k(\cdot,t)-u_{0k}\right\|_{L^1(\Re^n)}+\left\|u_{0k}-u_0\right\|_{L^1(\Re^n)}\\[0.5cm]\displaystyle
=2\left\|u_{0k}-u_0\right\|_{L^1(\Re^n)}+\left\|u_k(\cdot,t)-u_{0k}\right\|_{L^1(\Re^n)},\
(k=1,2,\cdots), \ \forall t>0.
\end{array}
\] By (\ref{4.4}) and (\ref{4.2}), we have
$$\displaystyle \lim\limits_{t\rightarrow 0^+}u(\cdot,t)=u_0\mbox{ in }L^1(\Re^n).$$
Combining this, (\ref{4.11}), (\ref{4.13}) with (\ref{4.10}),
we have proved the conclusion of Step 1. Moreover by the Steps 1 and 2
in the proof of of Theorem \ref{theo1}, the solution $u$ is unique, and satisfies
\begin{equation}\label{4.14}\begin{array}{l}\displaystyle
0\le V(y,\tau)\le F(y),\ \forall (y,\tau)\in Q,\\[0.3cm]\displaystyle \mbox{ and }
\left\|u(\cdot,t)\right\|_{L^\infty(\Re^n)}\le
Ct^{-\frac{1}{\beta}}\left\|u_0\right\|^{\frac{2}{\beta}}_{L^1(\Re^n)},
\ \forall t>0\end{array}\end{equation} by (\ref{4.6}),
(\ref{4.7}), (\ref{4.2}) and the  (\ref{3.7}) in Lemma 3.2

\underline{\bf Step 2.} Suppose $u_0\in L^1(\Re^n),\ u_0\ge 0$ in
$\Re^n.$ We will prove that (\ref{1.1}) has a nonnegative
$L^1$-regular solution $u$ satisfying $u(x,0)=u_0(x)$ in $\Re^n.$

Let $\chi_{B_k}$ be the characteristic function of the ball
$B_k(0)$ and set $$\displaystyle
u_{0k}(x)=\min\left\{k,u_0(x)\right\}\chi_{B_k}(x).$$ Then
\begin{equation}\displaystyle \label{4.15} 0\le u_{0k}\le
u_{0k'}\le u_0\mbox{ in }\Re^n \ \ (k\le k' , \ \  k, k' =1, 2,
\cdots,) u_0 =\lim\limits_{k\rightarrow\infty}u_{0k}\mbox{ in
}L^1(\Re^n).
\end{equation}
For each $k>0$, there exists a unique $L^1$-regular solution $u_k$
of (\ref{1.1}) such that $u_k(x,0)=u_{0k}$ in $\Re^n$ due to the
result of Step 1. In particular, we have
\begin{equation}\displaystyle \label{4.16}\int\limits_{\Re^n}u_k(x,t)dx=\int\limits_{\Re^n}u_{0k}(x)dx\le\left\|u_0\right\|_{L^1(\Re^n)},\ \forall t>0,\ \forall k=1,2,\cdots,\end{equation}
\begin{equation}\displaystyle \label{4.17}
\lim\limits_{t\rightarrow 0^+}u_{k}(\cdot,t)=u_{0k}\mbox{ in }L^1(\Re^n)
\ k=1,2,\cdots
\end{equation}
and
\begin{equation}\displaystyle \label{4.18}
0\le u_k(x,t)\le u_{k'}(x,t),\
\forall(x,t)\in\Re^n\times(0,\infty),\ \forall k\le
k'\end{equation} by the comparison principle as in the {\sl Steps}
1 and 2 of the proof of Theorem \ref{theo1}. Furthermore,
(\ref{4.13}) and (\ref{4.15}) imply
\begin{equation}\displaystyle \label{4.19}
\left\|u_k(\cdot,t)\right\|_{L^\infty(\Re^n)}\le
Ct^{-\frac{1}{\beta}}\left\|u_{0k}\right\|^{\frac{2}{\beta}}_{L^1(\Re^n)}
\le
Ct^{-\frac{1}{\beta}}\left\|u_{0}\right\|^{\frac{2}{\beta}}_{L^1(\Re^n)},
\ \forall t>0.\end{equation}

Now for any $t_0>0$, we use (\ref{4.19}) and Lemma 3.1 in [4] (or
more generally, Theorem 5.1 in [6]) to see that $\{u_k\}$ is
equi-continuous. Thus, (\ref{4.15})--(\ref{4.19}) imply that there
is a $u\in C(Q)\bigcap L^\infty(\Re^n\times[t_0,\infty))$ such
that $$\displaystyle \lim\limits_{k\rightarrow\infty}u_k=u\mbox{
in }Q,\
\lim\limits_{k\rightarrow\infty}u_k(\cdot,t)=u(\cdot,t)\mbox{ in
}L^1(\Re^n), \ \forall t>0,$$
\begin{equation}\displaystyle \label{4.20}
\left\|u(\cdot,t)\right\|_{L^\infty(\Re^n)} \le
Ct^{-\frac{1}{\beta}}\left\|u_{0}\right\|^{\frac{2}{\beta}}
_{L^1(\Re^n)},\ \
\int\limits_{\Re^n}u(x,t)dx=\int\limits_{\Re^n}u_0(x)dx,\ \forall
t>0\end{equation} and
\begin{equation}\displaystyle \label{4.21}
0\le u_k\le u\mbox{ in }Q,\ k=1,2,\cdots.\end{equation}
Obviously, $u$ is a solution of (\ref{1.1}). To complete the proof, we want only
 to prove $u\in C([0,\infty),L^1(\Re^n)).$ For this purpose, we fix
a $t_0>0.$ For any given $\varepsilon>0$ we use (\ref{4.15}) to see that
there exists a $k_0$ such that
\begin{equation}\displaystyle \label{4.22}
\int\limits_{\Re^n}(u_0(x)-u_{0k_0}(x)dx=\left\|u_0-u_{0k_0}\right\|\le
\frac{\varepsilon}{4}.\end{equation} Since $u_{k_0}$ is continuous
and $V_{k_0}$, related to $u_{k_0}$ as in Lemma 2.2, satisfies
(\ref{4.14}) for some $F\in L^1(\Re^n)$, we can find a $\delta>0$
such that
\begin{equation}\displaystyle \label{4.23}
\int\limits_{\Re^n}\left|u_{k_0}(\cdot,t)-u_{k_0}(\cdot,t_0)\right|dx\le
\frac{\varepsilon}{2},\ \forall t\in(t_0-\delta
,t_0+\delta).\end{equation} This, together with (\ref{4.20}),
(\ref{4.21}), (\ref{4.16}) and (\ref{4.22})
\[\begin{array}{l}
\int\limits_{\Re^n}\left|u(\cdot,t)-u(\cdot,t_0)\right|dx\\[0.5cm]\displaystyle
\le\int\limits_{\Re^n}\left|u(\cdot,t)-u_{k_0}(\cdot,t)\right|dx+
\int\limits_{\Re^n}\left|u_{k_0}(\cdot,t)-u_{k_0}(\cdot,t_0)\right|dx+
\int\limits_{\Re^n}\left|u_{k_0}(\cdot,t_0)-u(\cdot,t_0)\right|dx\\[0.5cm]\displaystyle
=2\int\limits_{\Re^n}\left(u_{0}(x)-u_{0k}(x)\right)dx+\int\limits_{\Re^n}
\left|u_{k_0}(\cdot,t)-u_{k_0}(\cdot,t_0)\right|dx\\[0.5cm]\displaystyle
\le \varepsilon,\ \ \forall t\in(t_0-\delta ,t_0+\delta),
\end{array}
\]
which shows that $u\in C((0,\infty),L^1(\Re^n))$. Repeating the
same argument for $t_0=0,$ with (\ref{4.17}) instead of
(\ref{4.23}), we see that $\lim\limits_{t\rightarrow
0^+}u(\cdot,t)=u_0$ in $L^1(\Re^n).$ This proves Theorem
\ref{theo2}.

\newpage
\section*{\centerline{REFERENCES}}
\begin{enumerate}
\item{}  D.G. Aronson, " The porous medium equation, In
{\em Nonlinear Diffusion Problems", Lecture Notes in Mathematics
1224,} (Edited by  Fasano A. and Primicerio M.), pp.1--46, Spinger
Verlag, New York, {\bf 1986}.

\item{}  E. DiBenedetto,  {\em Degenerate parabolic
equations}, Springer-Verlag, New York, {\bf 1993}.

\item{} J. Bear, {\em Dynamics of fluids in porous media,}
American Elsevier, New York, {\bf 1972}.

\item{} B. Song,  Anisotropic diffusions with singular
advections and absorptions, Part I, Existence, {\em Appl. Math.
Lett.},,
 {\bf 14} (2001), 811--816.

\item{} B. Song,  Anisotropic diffusions with singular
advections and absorptions, Part II, Uniqueness, {\em Appl. Math.
Lett.},{\bf 14} (2001), 817--823.

\item{}B. Song, On continuity study for some classes of
functions from degenerate parabolic equations and its
applications, to appear in {\em Northeastern Math. J.}.

\item{}H. Jian and B. Song, Fundamental solution of the
Anisotropic porous medium equation, to appear in {\em  Acta Math.
Sinica }.

\item{} Ph.  Benilan and  M. G.  Crandall,  The continuous
dependence on $\phi$ of solutions of $u_t-\Delta\phi(u)=0$, {\em
Indiana Univ. Math. J.} , {\bf 30 } (1981), 161-177.

\item{}  N. Alikakos  and  R. Rostamian, On the
uniformization of the solutions of the porous medium equation in
$\Re^n$, {\em Israel J. Math.},{\bf 47} (1984), 270--290.

\end{enumerate}

\end{document}